\documentclass[12pt,reqno]{amsart}
\usepackage{fullpage}
\usepackage{amssymb}
\usepackage{amsthm}
\usepackage{amsmath}
\usepackage[usenames]{color}
\usepackage{graphicx}
\usepackage[colorlinks=true,
linkcolor=webgreen,
filecolor=webbrown,
citecolor=webgreen]{hyperref}

\definecolor{webgreen}{rgb}{0,.5,0}
\definecolor{webbrown}{rgb}{.6,0,0}
\definecolor{red}{rgb}{1,0,0}

\usepackage{color}
\newtheorem{theorem}{Theorem}[section]

\newtheorem{lemma} [theorem]{Lemma}

\newtheorem{remark}[theorem]{Remark}

\begin{document}

\title{On the nonnegative integer solutions to the equation $F_n^2+F_m=2^{a}$}

\author{Sagar Mandal$^{*}$}
\address{Department of Mathematics\\Indian Institute of Technology Ropar, Punjab, India.}
\email{sagar.25maz0008@iitrpr.ac.in, sagarmandal31415@gmail.com}
\thanks{$^{*}$The author is financially supported by the University Grants Commission (UGC), Government of India, through the award of the Junior Research Fellowship (JRF) with ref. No. 241620111598.}

\maketitle
\let\thefootnote\relax
\footnotetext{\it MSC2020: 11B39, 11D61, 11J86} 
\maketitle
\footnotetext{\it Keywords: Diophantine equation, Fibonacci numbers, linear forms in logarithms, reduction method}
\let\thefootnote\relax

    \begin{abstract}
The Diophantine equation $F_n^x+F_m^y=2^a,$ where $F_k$ denotes the $k$-th Fibonacci number and $n,m,a\in \mathbb{Z}_{\geq0}$, has attracted attention in recent years. The cases $(x,y)=(1,1)$ and $(x,y)=(2,2)$ were completely resolved by Bravo and Luca, and by Jeong and Park, respectively. In this note, we investigate the mixed exponent case $(x,y)=(2,1)$ and we completely determine all nonnegative integer solutions.
\end{abstract}

\section{\protect\bigskip Introduction}
The Fibonacci sequence $\{F_n:n\geq0\}$ (\href{https://oeis.org/A000045}{\underline{A000045}}) is the binary recurrence sequence defined $F_0=0,F_1=1,$ and $ F_{n+2}=F_{n+1}+F_n$ for $n\geq0$. This is one of the most extensively studied sequences in number theory. In particular, Diophantine equations involving Fibonacci numbers and powers have attracted considerable attention over the years. A classical result of Bugeaud, Mignotte, and Siksek \cite{r4} completely determined the perfect powers appearing in the Fibonacci sequence. Later, Bravo and Luca \cite{r3} studied the equation $F_n+F_m=2^a$
and proved that it has only finitely many nonnegative integer solutions. Bravo and Bravo \cite{r1} extended this work to sums of three Fibonacci numbers, while Tiebekabe and Diouf \cite{r15,r16} considered sums of four and five Fibonacci numbers equal to powers of two. More generally, equations of the form $F_n\pm F_m=y^a$ have been investigated by several authors. In particular, Siar and Keskin \cite{r14} solved the equation $F_n-F_m=2^a,$
while Demirtürk Bitim and Keskin \cite{r6} studied the case $F_n-F_m=3^a.$
Further results for powers of five and for arbitrary perfect powers were obtained by Erduvan and Keskin \cite{r8}, Kebli, Kihel, Larone, Luca \cite{r10,r11}, and others. Recently, Jeong and Park \cite{r9} investigated the equation
$F_n^2+F_m^2=2^a$
and determined all its nonnegative integer solutions. Motivated by these developments and by the growing interest in exponential Diophantine equations involving Fibonacci numbers, we study the equation
\begin{align}\label{eq1}
    F_n^2+F_m=2^a.
\end{align}
\begin{remark}
For the Diophantine equations
$$
F_n+F_m=2^a \quad \text{and} \quad F_n^2+F_m^2=2^a,
$$
the equations are symmetric in the variables $n$ and $m$. Therefore, without loss of generality, one may assume that $n\ge m$. In contrast, the equation
$$
F_n^2+F_m=2^a
$$
is not symmetric, since the roles of $n$ and $m$ are fundamentally different. Consequently, we have to treat the two cases $n\ge m$ and $m>n$ separately.
\end{remark}
The OEIS sequences \href{https://oeis.org/A007598}{\underline{A007598}} and \href{https://oeis.org/A000079}{\underline{A000079}} correspond to the sequences $\{F_n^2:n\geq0\}$ and $\{2^n:n\geq0\}$, respectively.
Our main theorem  is stated below. 

\begin{theorem}\label{Thm1}
The Diophantine equation \eqref{eq1} has the following solutions in nonnegative integers $n,m,a$:
\begin{align*}
   \left( n,m,a\right) \in \big\{ &(1,0,0),(1,1,1),(2,0,0),(2,1,1),(2,2,1),(3,0,2),(6,0,6), (0,1,0), (0,2,0),
\\ &(0,3,1), (0,6,3), (1,4,2),(4,10,6), (1, 2, 1), (2, 4, 2) \big\}. 
\end{align*}

\end{theorem}

\section{Preliminaries}
We start this section by recalling some basic notions from algebraic number theory.

Let $\eta $ be an algebraic number of degree $d$ with minimal polynomial 
$$a_{0}x^{d}+a_{1}x^{d-1}+\ldots+a_{d}=a_{0}\prod\limits_{i=1}^{d}\left(x-\eta
^{(i)}\right), $$
where the $a_{i}$'s are relatively prime integers with $a_{0}>0$ and $\eta
^{(i)}$'s are conjugates of $\eta $. Then 
\begin{align*}
  h(\eta )=\frac{1}{d}\left( \log a_{0}+\sum\limits_{i=1}^{d}\log \left( \max
\left\{ |\eta ^{(i)}|,1\right\} \right) \right)  
\end{align*}
is called the logarithmic height of $\eta .$ In particular, if $\eta =a/b$
is a rational number with $\gcd (a,b)=1$ and $b>0,$ then $h(\eta )=\log
\left( \max \left\{ |a|,b\right\} \right)$. The following are some well known properties of logarithmic height:
\begin{enumerate}
    \item[\upshape(1)] $h(\eta \pm \gamma )\leq h(\eta )+h(\gamma )+\log 2$,
    \item[\upshape(2)] $h(\eta \gamma ^{\pm 1})\leq h(\eta )+h(\gamma )$,
    \item[\upshape(3)]$h(\eta ^{s})\leq |s|h(\eta ),~s\in \mathbb{Z}$. 
\end{enumerate}
In order to prove Theorem \ref{Thm1}, we need the following theorem due to Matveev. 
\begin{theorem}[\cite{r13}]\label{Thm2}
Assume that $\gamma _{1},\gamma _{2},...,\gamma _{r}$ are
positive real algebraic numbers in a real algebraic number field $\mathbb{K}$
of degree $D$, $b_{1},b_{2},\ldots,b_{r}\in \mathbb{Z}$, and 
$$\Lambda :=\gamma _{1}^{b_{1}}...\gamma _{r}^{b_{r}}-1 $$%
is not zero. Then 
$$|\Lambda |>\exp \left( -1.4\cdot 30^{r+3}\cdot r^{4.5}\cdot D^{2}(1+\log
D)(1+\log B)A_{1}A_{2}...A_{r}\right) , 
$$%
where 
$$
B\geq \max \left\{ |b_{1}|,...,|b_{r}|\right\} , 
$$%
and 
$$A_{i}\geq \max \left\{ Dh(\gamma _{i}),|\log \gamma _{i}|,0.16\right\}\quad \text{for all } i=1,...,r. $$

\end{theorem}
Next, we present the following lemma which is an immediate variation of a result due to Dujella and Peth{\H{o}} \cite{r7}. In the following lemma, the function $||\cdot ||$ denotes the
distance from $x$ to the nearest integer, that is, $||x||=\min \left\{
|x-n|:n\in\mathbb{Z}\right\} $ for a real number $x.$

\begin{lemma}\label{Lem1}Let $M$ be a positive integer, let $p/q$ be a convergent of the
continued fraction expansion of the irrational number $\gamma $ such that $q>6M,$ and
let $A,B,\mu $ be some real numbers with $A>0$ and $B>1.$ Let $\epsilon
:=||\mu q||-M||\gamma q||.$ If $\epsilon >0,$ then there exists no solution
to the inequality 
$$
0<|u\gamma -v+\mu |<AB^{-w}, 
$$
in positive integers $u,v,$ and $w$ with 
$$
u\leq M\text{ and }w\geq \frac{\log (Aq/\epsilon )}{\log B}. 
$$
\end{lemma}
The following lemma will play a pivotal role in subsequent proofs.
\begin{lemma}[{\cite{rnew}}]\label{lemma2.3}
For all $n \ge 1$, we have 
$$
\nu_2(F_n)=
\begin{cases}
0, & \text{if } n \equiv 1,2 \pmod{3},\\[4pt]
1, & \text{if } n \equiv 3 \pmod{6},\\[4pt]
3, & \text{if } n \equiv 6 \pmod{12},\\[4pt]
\nu_2(n)+2, & \text{if } n \equiv 0 \pmod{12}.
\end{cases}
$$
\end{lemma}
It is well known that 
\begin{align}\label{eq2}
F_{n}=\dfrac{\phi ^{n}-\psi ^{n}}{\sqrt{5}}   \quad \text{for  } n\geq0, 
\end{align}
where $\phi =\dfrac{1+\sqrt{5}}{2}$ and $\psi =\dfrac{1-\sqrt{5}}{2},$ which are the roots of the characteristic equation $x^{2}-x-1=0$. The following inequality can be proved by induction:
\begin{align}\label{eq3}
    \phi ^{n-2}\leq F_{n}\leq \phi ^{n-1}  \quad \text{for  } n\geq1.
\end{align}
\section{Proof of Theorem \ref{Thm1}}
First, observe that equation (\ref{eq1}) has no solution when $n=m=0$. We divide the proof into the following three parts.
Part-1: Let $2n>m$. We consider the following two cases:\\
Case-1: For $m=0$ we have $F_n^2=2^{a}$, from \cite{r9}, we get the following solutions
$$\left( n,m,a\right) \in \big\{(1,0,0), (2,0,0), (3,0,2), (6,0,6)\big\}.$$
Case-2: Let $2n> m\geq 1$. For $1\leq n\leq 5$, (\ref{eq1}) has the following solutions
$$\left( n,m,a\right) \in \big\{(1,1,1),(2,1,1),(2,2,1)\big\}.$$

First, we establish an upper bound for $n\geq 6$ using Theorem \ref{Thm2}. 

From (\ref{eq3}), we have
$2^{a}=F^2_{n}+F_{m}\leq \phi^{2n}+\phi^{m}<2^{2n}+2^{2n}=2^{2n+1}.$
It follows that $a\leq 2n$. Now, using (\ref{eq2}), $F^2_{n}+F_{m}=2^{a}$ can be written as
\begin{align*}\label{eq4}
\big(\frac{\phi^n-\psi^n}{\sqrt{5}}\big)^2+\frac{\phi^m-\psi^m}{\sqrt{5}}=2^{a}.\tag{4}    
\end{align*}

From the above equation and using the fact that $|\psi|<1$, we obtain
\begin{align*}\label{eq5}
    \big|\frac{\phi^{2n}}{5}-2^{a}\big|=\big|\frac{2\phi^n\psi^n}{5}-\frac{\psi^{2n}}{5}+\frac{\psi^{m}}{\sqrt{5}}-\frac{\phi^{m}}{\sqrt{5}}\big|&< \frac{2\phi^{n+m/2}}{5}+\frac{\phi^{n+m/2}}{5}+\frac{\phi^{n+m/2}}{\sqrt{5}}+\frac{\phi^{n+m/2}}{\sqrt{5}}\\
    &<\frac{8\phi^{n+m/2}}{5}.\tag{5}
\end{align*}
After multiplying both sides of the above inequality by $\dfrac{5}{\phi^{2n}}$, we get
\begin{align*}\label{eq6}
    \big|1-2^{a}\cdot \frac{5}{\phi^{2n}}\big|<\frac{8}{\phi^{n-m/2}}.\tag{6}
\end{align*}
Now, let us apply Theorem \ref{Thm2} with $(\gamma _{1},b_1)=(2,a),~(\gamma _{2},b_2)=(\sqrt{5} ,2),~(\gamma_{3},b_3)=(\phi,-2n)$, and 
$$\Lambda_1 = 2^{a}\cdot \frac{5}{\phi^{2n}} - 1.$$
Observe that $\Lambda_1\neq0$, as if $\Lambda_1=0$, we have $\phi^{2n}=2^{a}\cdot 5$, then conjugating the equation in $\mathbb{Q}(\sqrt{5})$, we get $\psi^{2n}=2^{a}\cdot 5$, which is absurd because $|\psi|<1$. Also, note that $D=[\mathbb{Q}(\sqrt5):\mathbb{Q}]=2$. Set
$$B=2n\geq\max\{|b_1|,|b_2|,|b_3|\},\quad \text{and} \quad A_1=1.39\geq \max\{D h(\gamma_1),\,|\log \gamma_1|,\,0.16\}=2\log2.$$
Since the minimal polynomial of $\gamma_2=\sqrt5$ is $x^2-5$ and $h(\gamma_2)=\frac12\log 5$, we set
$$A_2=1.61\geq \max\{D h(\gamma_2),\,|\log \gamma_2|,\,0.16\}=\log 5.$$
Since the minimal polynomial of $\gamma_3=\phi$ is $x^2-x-1$ and $h(\gamma_3)=\frac{1}{2}\log\phi $, we choose 
$$A_3=0.49\geq \max\{D h(\gamma_3),\,|\log \gamma_3|,\,0.16\}=\log \phi.$$

Then, by Theorem \ref{Thm2} and (\ref{eq6}), we obtain 
\begin{align*}
    \exp(-1.4\cdot 30^{6}\cdot 3^{4.5}\cdot D^2(1+\log D)\cdot A_1\cdot A_2\cdot A_3\cdot (1+\log 2n))<|\Lambda_1|<\frac{8}{\phi^{n-m/2}}.
\end{align*}
It follows that
\begin{align*}
-1.4\cdot 30^{6}\cdot 3^{4.5}\cdot D^2(1+\log D)\cdot A_1\cdot A_2\cdot A_3\cdot(1+\log 2n)<\log8-(n-m/2)\log \phi,
\end{align*}
from which we get
\begin{align*}\label{eq7}
0<(n-m/2)\log \phi&<\log8+1.4\cdot 30^{6}\cdot 3^{4.5}\cdot D^2(1+\log D)\cdot A_1\cdot A_2\cdot A_3\cdot(1+\log 2n)\\
&<11\cdot 10^{11}\cdot (1+\log 2n).\tag{7}
\end{align*}
Further, using (\ref{eq4}), we can write
\begin{align*}\label{eq8}
    \big|\frac{\phi^{2n}}{5}+\frac{\phi^m}{\sqrt{5}}-2^{a}\big|=\big|\frac{2\phi^n\psi^n}{5}-\frac{\psi^{2n}}{5}+\frac{\psi^{m}}{\sqrt{5}}\big|<\frac{2\phi^{n}}{5}+\frac{\phi^{n}}{5}+\frac{\phi^{n}}{\sqrt{5}}<\frac{11\phi^{n}}{10}.\tag{8}
\end{align*}
After multiplying both sides of the above inequality by $\dfrac{5}{\phi^{2n}+\sqrt{5}\phi^m}$, we get
\begin{align*}\label{eq9}
    \big|1-2^{a}\cdot \frac{5}{\phi^{2n}+\sqrt{5}\phi^m}\big|<\frac{6}{\phi^{n}}.\tag{9}
\end{align*}
In Theorem \ref{Thm2}, set $(\gamma _{1},b_1)=(2,a),~(\gamma _{2},b_2)=(\sqrt{5} ,2),~(\gamma_{3},b_3)=(\phi,-m),~(\gamma_4,b_4)=(\phi^{2n-m}+\sqrt{5},-1)$, and 
$$\Lambda_2 = 2^{a}\cdot \frac{5}{\phi^{2n}+\sqrt{5}\phi^m} - 1.$$
Observe that $\Lambda_2\neq0$, as if $\Lambda_2=0$, we have $\phi^{2n}+\sqrt{5}\phi^m=2^{a}\cdot 5$, then conjugating the equation in $\mathbb{Q}(\sqrt{5})$, we get $\psi^{2n}-\sqrt{5}\psi^{m}=2^{a}\cdot 5$, which is absurd because 
$$2^{a}\cdot 5\leq |\psi|^{2n}+\sqrt{5}|\psi|^m<1+\sqrt{5}<5.$$ 
Now, choose
$$B=2n\geq \max\{|b_1|,|b_2|,|b_3|,|b_4|\},\quad \quad A_1=1.39\geq \max\{D h(\gamma_1),\,|\log \gamma_1|,\,0.16\}=2\log2,$$
$$A_2=1.61\geq \max\{D h(\gamma_2),\,|\log \gamma_2|,\,0.16\}=\log 5,\quad A_3=0.49\geq \max\{D h(\gamma_3),\,|\log \gamma_3|,\,0.16\}=\log \phi.$$
Using properties of  logarithmic height, we have
$$h(\gamma_4)\leq h(\sqrt{5})+\log 2+(2n-m)h(\phi)=\frac{1}{2}\log5+\log 2+(n-m/2)\log \phi<1.5+(n-m/2)\log \phi.$$
Therefore, we choose
\begin{align*}
A_4=3+(2n-m)\log \phi\geq \max\{D h(\gamma_4),\,|\log \gamma_4|,\,0.16\}. 
\end{align*}
From Theorem \ref{Thm2} and (\ref{eq9}), we obtain
\begin{align*}
    \exp(-1.4\times 30^{\,7}\cdot 4^{4.5}\cdot D^2(1+\log D)\cdot A_1\cdot A_2\cdot A_3\cdot A_4\cdot (1+\log 2n))<|\Lambda_2|<\frac{6}{\phi^{n}}.
\end{align*}
It follows that
\begin{align*}
-1.4\times 30^{\,7}\cdot 4^{4.5}\cdot D^2(1+\log D)\cdot A_1\cdot A_2\cdot A_3\cdot A_4\cdot (1+\log 2n)<\log6-n\log \phi,
\end{align*}
from which we get
\begin{align*}
n\log \phi&<\log6+1.4\times 30^{\,7}\cdot 4^{4.5}\cdot D^2(1+\log D)\cdot A_1\cdot A_2\cdot A_3\cdot A_4\cdot (1+\log 2n)\\
&<12\cdot 10^{13}\cdot \big(3+(2n-m)\log\phi\big)(1+\log 2n)\tag{10}\label{eq10}\\ 
&<12\cdot 10^{13}\cdot \big(3+22\cdot 10^{11}(1+\log 2n)\big)(1+\log 2n),\tag{11}\label{eq11}
\end{align*}
the last inequality follows from (\ref{eq7}). The above inequality does not hold for all $n>2.85\cdot 10^{30}$. Thus, we must have $n\leq 2.85\cdot 10^{30}$. We now need to reduce the upper bound obtained for $n$. For this, we use Lemma \ref{Lem1}.\\
Let $z_1=a\log 2-2n\log \phi+\log 5$.
Observe that $z_1>0$, if not, then $a\log2+\log5\leq 2n\log\phi$, which implies that $2^{a}\cdot 5\leq \phi^{2n}$, which is impossible as $2^{a}\leq \dfrac{\phi^{2n}}{5}\leq F_n^2+\dfrac{2}{5}-\dfrac{\psi^{2n}}{5}<F_n^2+F_m=2^{a}$. From (\ref{eq6}), we have
$$0<z_1\leq e^{z_1}-1<\frac{8}{\phi^{n-m/2}},$$
from which we can write
\begin{align*}
    a\log 2-2n\log \phi+\log 5<\frac{8}{\phi^{n-m/2}},
\end{align*}
that is
\begin{align*}\label{eq12}
a\frac{\log 2}{\log \phi}-2n+\frac{\log 5}{\log \phi}<\frac{8\phi^{m/2-n}}{\log \phi}<17\phi ^{m/2-n}.\tag{12}
\end{align*}
Now, let us apply Lemma \ref{Lem1} with
\begin{align*}
    u=a,\quad \gamma=\frac{\log2}{\log \phi},\quad v=2n,\quad \mu=\frac{\log5}{\log \phi},\quad A=17,\quad B=\sqrt{\phi},\quad w=2n-m.
\end{align*}
Since $n\leq 2.85\cdot 10^{30}$ and $u=a\leq 2n<5.7\cdot 10^{30}$, we set $M=5.7\cdot 10^{30}$. Choosing  
$$\frac{p}{q}=\frac{p_{70}}{q_{70}}=\frac{228666343422267608843910896109913}{158749759840390984049158390593424},$$
we get $q_{70}>6M$ and
$$\epsilon=||\mu q||-M||\gamma q||=0.0182\ldots>0.$$
Therefore, by Lemma \ref{Lem1}, we have
\begin{align*}\label{eq13}
   2n-m<\frac{\log(Aq/\epsilon)}{\log B}<337.\tag{13}
\end{align*}
Now, using (\ref{eq10}) and (\ref{eq13}), we get
\begin{align*}
n\log \phi<12\cdot 10^{13}\cdot \big(3+(2n-m)\log\phi\big)(1+\log 2n)\leq12\cdot 10^{13}\cdot \big(3+336\cdot\log \phi\big)(1+\log 2n).
\end{align*}
The above inequality implies that $1\leq n<1.8\cdot 10^{18}$.

Again, let
$$z_2=a\log 2-m\log\phi+\log \frac{5}{\phi^{2n-m}+\sqrt{5}}.$$
Note that $z_2\neq 0$. If $z_2>0$, then from (\ref{eq9}), we have
\begin{align*}\label{eq14}
0<z_2\leq e^{z_2}-1<\frac{6}{\phi^{n}}. \tag{14}   
\end{align*}
If $z_2<0$, then from (\ref{eq9}), we have
\begin{align*}\label{eq15}
    |e^{z_2}-1|\leq \frac{6}{\phi^n}<\frac{1}{2},\quad \text{for }n\geq6.\tag{15}
\end{align*}
Recall that for any nonzero real number $x<0$, if $|e^{x}-1| <\dfrac{1}{2}$, then $
|x| < 2|e^{x}-1|$. Therefore, for $n\geq6$, we obtain from (\ref{eq14}), (\ref{eq15}) that
\begin{align*}\label{eq16}
    0<|z_2|<\frac{12}{\phi^n}.\tag{16}
\end{align*}
For $x\geq0$, set $f(x)=\dfrac{5}{\phi^{x}+\sqrt{5}}$. Then, from (\ref{eq16}), we can write
\begin{align*}
    |a\log 2-m\log\phi+\log f(2n-m)|<\frac{12}{\phi^{n}},
\end{align*}
multiplying both sides of the above inequality by $\dfrac{1}{\log\phi}$, we get
\begin{align*}
    \bigg|a\frac{\log 2}{\log \phi}-m+\frac{\log f(2n-m)}{\log \phi}\bigg|<\frac{12\phi^{-n}}{\log \phi}<25\phi ^{-n}.
\end{align*}
In Lemma \ref{Lem1}, we set
\begin{align*}
    u=a,\quad \gamma=\frac{\log2}{\log \phi},\quad v=m,\quad \mu=\mu_{2n-m}=\frac{\log f(2n-m) }{\log \phi},\quad A=25,\quad B=\phi,\quad w=n.
\end{align*}
Since $n<1.8\cdot 10^{18}$ and $u=a\leq 2n<3.6\cdot 10^{18}$, we set $M=3.6\cdot 10^{18}$. Then choosing 

$$\frac{p}{q}=\frac{p_{70}}{q_{70}}=\frac{228666343422267608843910896109913}{158749759840390984049158390593424},$$
we get $q_{70}>6M$, and
$$\epsilon_{2n-m}=||\mu_{2n-m} q||-M||\gamma q||>0,\quad \text{for all }2n-m=1,2,\ldots,336.$$
Therefore, applying Lemma \ref{Lem1}, we can see that for each $2n-m=1,2,\ldots,336$, the following inequality holds
\begin{align*}
    n<\frac{\log(Aq/\epsilon_{2n-m})}{\log B}\leq \max_{\epsilon_{2n-m}} \frac{\log(Aq/\epsilon_{2n-m})}{\log B}<179.
\end{align*}
After running a brute force search using Python code for $n$ and $m$ with $1\leq m<2n\leq 356$, $n\geq 6$ satisfying (\ref{eq1}) for some nonnegative integer $a$, we obtain no solutions.

Part-2: Let $m>2n$. We consider the following two cases:

Case-1: Let $n=0$. Then we have $F_m=2^{a}$. If $a=0$, then $F_m=1$ for $m=1,2$. If $a\geq1$ then from Lemma \ref{lemma2.3}, we have $3\mid m$. Now for such $m$, we must have either $F_m=2$ or $F_m=2^{\nu_2(m)+2}$. Now $F_m=2$ holds for $m=3$. For $m\geq6$, if $F_m=2^{\nu_2(m)+2}$ then $F_m\leq4m$ as $\nu_2(m)\leq \log_2m$ for $m\geq3$, but $F_m\leq 4m$ is true for $m=6,9$. Hence, if $n=0$ and $m\geq1$, (\ref{eq1}) has the following solutions
$$\left( n,m,a\right) \in \big\{(0,1,0), (0,2,0), (0,3,1), (0,6,3)\big\}.$$

Case-2: Let $n\geq1$. For $2<m\leq 11$, (\ref{eq1}) has the following solutions

$$\left( n,m,a\right) \in \big\{(1,4,2),(4,10,6)\big\}.$$
Therefore, we consider $m\geq12$. From (\ref{eq3}), we have
$2^{a}=F^2_{n}+F_{m}\leq \phi^{2n}+\phi^{m}\leq \phi^m+\phi^m<2^{m+1}.$
It follows that $a\leq m$. Now, using (\ref{eq4}), we have
\begin{align*}
    \big|\frac{\phi^{m}}{\sqrt{5}}-2^{a}\big|=\big|\frac{2\phi^n\psi^n}{5}-\frac{\phi^{2n}}{5}-\frac{\psi^{2n}}{5}+\frac{\psi^{m}}{\sqrt{5}}\big|<\phi^{2n}.
\end{align*}
After multiplying both sides of the above inequality by $\dfrac{\sqrt{5}}{\phi^{m}}$, we get
\begin{align*}\label{eq17}
    \big|1-2^{a}\cdot \frac{\sqrt{5}}{\phi^{m}}\big|<\frac{3}{\phi^{m-2n}}.\tag{17}
\end{align*}

Now, let us apply Theorem \ref{Thm2} with $(\gamma _{1},b_1)=(2,a),~(\gamma _{2},b_2)=(\sqrt{5} ,1),~(\gamma_{3},b_3)=(\phi,-m)$, and 
$$\Lambda_1 = 2^{a}\cdot \frac{\sqrt{5}}{\phi^{m}} - 1.$$
Observe that $\Lambda_1\neq0$, as if $\Lambda_1=0$, we have $\phi^{m}=2^{a}\cdot \sqrt{5}$, then conjugating the equation in $\mathbb{Q}(\sqrt{5})$, we get $\psi^{m}=-2^{a}\cdot \sqrt{5}$, which is absurd because $2^{a}\cdot \sqrt{5}=|\psi|^m<1$. Also, note that $D=[\mathbb{Q}(\sqrt5):\mathbb{Q}]=2$. Set
$$B=m\geq\max\{|b_1|,|b_2|,|b_3|\},\quad A_1=1.39\geq \max\{D h(\gamma_1),\,|\log \gamma_1|,\,0.16\}=2\log2.$$
$$A_2=1.61\geq \max\{D h(\gamma_2),\,|\log \gamma_2|,\,0.16\}=\log 5,\quad A_3=0.49\geq \max\{D h(\gamma_3),\,|\log \gamma_3|,\,0.16\}=\log \phi.$$
Then, by Theorem \ref{Thm2} and (\ref{eq17}), we obtain 
\begin{align*}
    \exp(-1.4\cdot 30^{6}\cdot 3^{4.5}\cdot D^2(1+\log D)\cdot A_1\cdot A_2\cdot A_3\cdot (1+\log m))<|\Lambda_1|<\frac{3}{\phi^{m-2n}}.
\end{align*}
It follows that
\begin{align*}
-1.4\cdot 30^{6}\cdot 3^{4.5}\cdot D^2(1+\log D)\cdot A_1\cdot A_2\cdot A_3\cdot(1+\log m)<\log3-(m-2n)\log \phi,
\end{align*}
from which we get
\begin{align*}\label{eq18}
0<(m-2n)\log \phi&<\log3+1.4\cdot 30^{6}\cdot 3^{4.5}\cdot D^2(1+\log D)\cdot A_1\cdot A_2\cdot A_3\cdot(1+\log m)\\
&<10.7\cdot 10^{11}\cdot (1+\log m).\tag{18}
\end{align*}

Further, using (\ref{eq4}), we can write
\begin{align*}
    \big|\frac{\phi^{2n}}{5}+\frac{\phi^m}{\sqrt{5}}-2^{a}\big|=\big|\frac{2\phi^n\psi^n}{5}-\frac{\psi^{2n}}{5}+\frac{\psi^{m}}{\sqrt{5}}\big|<\frac{2\phi^{m/2}}{5}+\frac{\phi^{m/2}}{5}+\frac{\phi^{m/2}}{\sqrt{5}}<\frac{11\phi^{m/2}}{10}.
\end{align*}
After multiplying both sides of the above inequality by $\dfrac{5}{\phi^{2n}+\sqrt{5}\phi^m}$, we get
\begin{align*}\label{eq19}
    \big|1-2^{a}\cdot \frac{5}{\phi^{2n}+\sqrt{5}\phi^m}\big|<\frac{6}{\phi^{m/2}}.\tag{19}
\end{align*}
In Theorem \ref{Thm2}, set $(\gamma _{1},b_1)=(2,a),~(\gamma _{2},b_2)=(\sqrt{5} ,2),~(\gamma_{3},b_3)=(\phi,-2n),~(\gamma_4,b_4)=(\sqrt{5}\phi^{m-2n}+1,-1)$, and 
$$\Lambda_2 = 2^{a}\cdot \frac{5}{\phi^{2n}+\sqrt{5}\phi^m} - 1.$$
As observed previously that $\Lambda_2\neq0$. Now, choose
$$B=m\geq \max\{|b_1|,|b_2|,|b_3|,|b_4|\},\quad \quad A_1=1.39\geq \max\{D h(\gamma_1),\,|\log \gamma_1|,\,0.16\}=2\log2,$$
$$A_2=1.61\geq \max\{D h(\gamma_2),\,|\log \gamma_2|,\,0.16\}=\log 5,\quad A_3=0.49\geq \max\{D h(\gamma_3),\,|\log \gamma_3|,\,0.16\}=\log \phi.$$
Using properties of  logarithmic height, we have
$$h(\gamma_4)\leq \frac{1}{2}\log5+\log 2+(m/2-n)\log \phi<1.5+(m/2-n)\log \phi.$$
Therefore, we choose
\begin{align*}
A_4=3+(m-2n)\log \phi\geq \max\{D h(\gamma_4),\,|\log \gamma_4|,\,0.16\}. 
\end{align*}
From Theorem \ref{Thm2} and (\ref{eq19}), we obtain
\begin{align*}
    \exp(-1.4\times 30^{\,7}\cdot 4^{4.5}\cdot D^2(1+\log D)\cdot A_1\cdot A_2\cdot A_3\cdot A_4\cdot (1+\log m))<|\Lambda_2|<\frac{6}{\phi^{m/2}}.
\end{align*}
It follows that
\begin{align*}
-1.4\times 30^{\,7}\cdot 4^{4.5}\cdot D^2(1+\log D)\cdot A_1\cdot A_2\cdot A_3\cdot A_4\cdot (1+\log m)<\log6-\frac{m}{2}\log \phi,
\end{align*}
from which we get
\begin{align*}
\frac{m}{2}\log \phi&<\log6+1.4\times 30^{\,7}\cdot 4^{4.5}\cdot D^2(1+\log D)\cdot A_1\cdot A_2\cdot A_3\cdot A_4\cdot (1+\log m)\\
&<12\cdot 10^{13}\cdot \big(3+(m-2n)\log\phi\big)(1+\log m)\tag{20}\label{eq20}\\ 
&<12\cdot 10^{13}\cdot \big(3+10.7\cdot 10^{11}(1+\log m)\big)(1+\log m),\tag{21}\label{eq21}
\end{align*}
the last inequality follows from (\ref{eq18}). From the above inequality, we must have $m\leq 2.7\cdot 10^{30}$. We now need to reduce the upper bound obtained for $m$. For this, as previously, we  use Lemma \ref{Lem1}.\\

Let $z_1=a\log 2-m\log \phi+\log \sqrt{5}$.
Observe that $z_1>0$, if not, then $a\log2+\log\sqrt{5}\leq m\log\phi$, which implies that $2^{a}\cdot \sqrt{5}\leq \phi^{m}$, which is impossible as 
$$F_n^2+\frac{\phi^m-\psi^m}{\sqrt{5}}=2^{a}\leq \dfrac{\phi^{m}}{\sqrt{5}},$$ 
which gives $F_n^2\leq \dfrac{\psi^m}{\sqrt{5}}$. From (\ref{eq17}), we have
$$0<z_1\leq e^{z_1}-1<\frac{3}{\phi^{m-2n}},$$
from which we can write
\begin{align*}
    a\log 2-m\log \phi+\log \sqrt{5}<\frac{3}{\phi^{m-2n}},
\end{align*}
that is
\begin{align*}
a\frac{\log 2}{\log \phi}-m+\frac{\log \sqrt{5}}{\log \phi}<7\phi ^{-(m-2n)}.
\end{align*}
Now, let us apply Lemma \ref{Lem1} with
\begin{align*}
    u=a,\quad \gamma=\frac{\log2}{\log \phi},\quad v=m,\quad \mu=\frac{\log\sqrt5}{\log \phi},\quad A=7,\quad B=\phi,\quad w=m-2n.
\end{align*}
Since $m\leq 2.7\cdot 10^{30}$ and $u=a\leq m\leq 2.7\cdot 10^{30}$, we set $M=2.7\cdot 10^{30}$. Choosing  
$$\frac{p}{q}=\frac{p_{70}}{q_{70}}=\frac{228666343422267608843910896109913}{158749759840390984049158390593424},$$
we get $q_{70}>6M$ and
$$\epsilon=||\mu q||-M||\gamma q||=0.0099\ldots>0.$$
Therefore, by Lemma \ref{Lem1}, we have
\begin{align*}\label{eq22}
    m-2n<\frac{\log(Aq/\epsilon)}{\log B}<168.\tag{22}
\end{align*}
Now, using (\ref{eq20}) and (\ref{eq22}), we get
\begin{align*}
\frac{m}{2}\log \phi&<\log6+1.4\times 30^{\,7}\cdot 4^{4.5}\cdot D^2(1+\log D)\cdot A_1\cdot A_2\cdot A_3\cdot A_4\cdot (1+\log m)\\
&<12\cdot 10^{13}\cdot \big(3+167\log\phi\big)(1+\log m),
\end{align*}
which is true only for $m<1.79\cdot 10^{18}$. Let us apply Lemma \ref{Lem1} one more time with
\begin{align*}
    u=a,\quad \gamma=\frac{\log2}{\log \phi},\quad v=m,\quad \mu=\frac{\log\sqrt5}{\log \phi},\quad A=7,\quad B=\phi,\quad w=m-2n.
\end{align*}
Since $m\leq 1.79\cdot 10^{18}$ and $u=a\leq m\leq 1.79\cdot 10^{18}$, we set $M=1.79\cdot 10^{18}$. Choosing  
$$\frac{p}{q}=\frac{p_{45}}{q_{45}}=\frac{78467056620067185965}{54475119544877440894},$$
we get $q_{45}>6M$ and
$$\epsilon=||\mu q||-M||\gamma q||>0.$$
Therefore, by Lemma \ref{Lem1}, we have
\begin{align*}\label{eq23}
    m-2n<\frac{\log(Aq/\epsilon)}{\log B}<102.\tag{23}
\end{align*}
Now, using (\ref{eq20}) and (\ref{eq23}), we get
\begin{align*}
\frac{m}{2}\log \phi&<\log6+1.4\times 30^{\,7}\cdot 4^{4.5}\cdot D^2(1+\log D)\cdot A_1\cdot A_2\cdot A_3\cdot A_4\cdot (1+\log m)\\
&<12\cdot 10^{13}\cdot \big(3+101\log\phi\big)(1+\log m),
\end{align*}
which gives $m<1.1\cdot 10^{18}$. 
Now, set
$$z_2=a\log 2-2n\log\phi+\log \frac{5}{\sqrt{5}\phi^{m-2n}+1}.$$
Note that $z_2\neq 0$. If $z_2>0$, then from (\ref{eq19}), we have
\begin{align*}\label{eq24}
0<z_2\leq e^{z_2}-1<\frac{6}{\phi^{m/2}}. \tag{24}   
\end{align*}
If $z_2<0$, then from (\ref{eq19}), we have
\begin{align*}\label{eq25}
    |e^{z_2}-1|\leq \frac{6}{\phi^{m/2}}<\frac{1}{2},\quad \text{for }m\geq12.\tag{25}
\end{align*}
Therefore, for $m\geq12$, we obtain from (\ref{eq24}), and (\ref{eq25}) that
\begin{align*}\label{eq26}
    0<|z_2|<\frac{12}{\phi^{m/2}}.\tag{26}
\end{align*}
For $x\geq0$, set $f(x)=\dfrac{5}{\sqrt{5}\phi^{x}+1}$. Then, from (\ref{eq26}), we can write
\begin{align*}
    |a\log 2-2n\log\phi+\log f(m-2n)|<\frac{12}{\phi^{m/2}},
\end{align*}
multiplying both sides of the above inequality by $\dfrac{1}{\log\phi}$, we get
\begin{align*}
    \bigg|a\frac{\log 2}{\log \phi}-2n+\frac{\log f(m-2n)}{\log \phi}\bigg|<\frac{12\phi^{-{m/2}}}{\log \phi}<25\phi ^{-m/2}.
\end{align*}
In Lemma \ref{Lem1}, we set
\begin{align*}
    u=a,\quad \gamma=\frac{\log2}{\log \phi},\quad v=2n,\quad \mu=\mu_{m-2n}=\frac{\log f(m-2n) }{\log \phi},\quad A=25,\quad B=\sqrt{\phi},\quad w=m.
\end{align*}
Since $m<1.1\cdot 10^{18}$ and $u=a\leq m<1.1\cdot 10^{18}$, we set $M=1.1\cdot 10^{18}$. Then choosing 

$$\frac{p}{q}=\frac{p_{50}}{q_{50}}=\frac{2697764232477112944722}{1872901003279144420473},$$
we get $q_{50}>6M$, and
$$\epsilon_{m-2n}=||\mu_{m-2n} q||-M||\gamma q||>0,\quad \text{for all }m-2n=1,2,\ldots,101.$$
Therefore, applying Lemma \ref{Lem1}, we can see that for each $m-2n=1,2,\ldots,101$, the following inequality holds
\begin{align*}
    m<\frac{\log(Aq/\epsilon_{m-2n})}{\log B}\leq \max_{\epsilon_{m-2n}} \frac{\log(Aq/\epsilon_{m-2n})}{\log B}<242.
\end{align*}
After running a brute force search using Python code for $n$ and $m$ with $2n<m<242$, $n\geq1$, and $m\geq12$ satisfying (\ref{eq1}) for some nonnegative integer $a$, we obtain no solutions.

Part-3: Let $2n=m$. Then (\ref{eq1}) becomes $F_n^2+F_{2n}=2^{a}$ which is equivalent to $2F_{n+1}F_{n}=2^{a}$ that is $F_{n+1}F_{n}=2^{a-1}$. If $a=1$, then $n=1$. For $n\geq2$, since $\gcd(F_n,F_{n+1})=1$, for $a\geq2$, we must have either $F_{n+1}=2^{a-1}$ with $F_n=1$ or $F_n=2^{a-1}$ with $F_{n+1}=1$.  Observe that $F_n=2^{a-1}$ with $F_{n+1}=1$ is impossible, and $F_{n+1}=2^{a-1}$ with $F_n=1$ holds for $a=2,n=2$. Therefore, for $2n=m$, the solutions are 
$$\left( n,m,a\right) \in \big\{(1,2,1),(2,4,2)\big\}.$$

This completes the proof.\qed

\section{Data Availability} 	
The author confirms that the manuscript has no associated data.

\section{Competing Interests}
The author confirms that he has no competing interest.


\begin{thebibliography}{99}
\bibitem{r1}
E. F. Bravo and J. J. Bravo,
\emph{Powers of two as sums of three Fibonacci numbers},
Lith. Math. J. \textbf{55} (2015), no. 3, 301--311.

\bibitem{r2}
J. J. Bravo, M. D{\'\i}az, and F. Luca,
\emph{Ratios of sums of two Fibonacci numbers equal to powers of 2},
Math. Commun. \textbf{25} (2020), no. 2, 185--199.

\bibitem{r3}
J. J. Bravo and F. Luca,
\emph{On the Diophantine equation $F_n+F_m=2^a$},
Quaest. Math. \textbf{39} (2016), no. 3, 391--400.

\bibitem{r4}
Y. Bugeaud, M. Mignotte, and S. Siksek,
\emph{Classical and modular approaches to exponential Diophantine equations. I. Fibonacci and Lucas perfect powers},
Ann. of Math. (2) \textbf{163} (2006), no. 3, 969--1018.
\bibitem{r5} K.C. Chim, V. Ziegler. On Diophantine equations involving sums of Fibonacci numbers and power of 2. Integers 2018; 18, Paper No. A99, 30 pp.


\bibitem{r6}
B. Demirt{\"u}rk Bitim and R. Keskin,
\emph{On solutions of the Diophantine equation $F_n-F_m=3^a$},
Proc. Indian Acad. Sci. Math. Sci. \textbf{129} (2019), no. 5, Paper No. 81, 10 pp.

\bibitem{r7}
A. Dujella and A. Peth{\H{o}},
\emph{A generalization of a theorem of Baker and Davenport},
Quart. J. Math. Oxford Ser. (2) \textbf{49} (1998), no. 195, 291--306.

\bibitem{r8}
F. Erduvan and R. Keskin,
\emph{Nonnegative integer solutions of the equation $F_n-F_m=5^a$},
Turkish J. Math. \textbf{43} (2019), no. 3, 1115--1123.
\bibitem{r9}
G. H. Jeong and P. S. Park,
\emph{On the Diophantine equation $F_n^2+F_m^2=2^a$},
Math. Commun. \textbf{30} (2025), no. 2, 171--177.
\bibitem{r10}
S. Kebli, O. Kihel, J. Larone, and F. Luca,
\emph{On the nonnegative integer solutions to the equation $F_n\pm F_m=y^a$},
J. Number Theory \textbf{220} (2021), 107--127.

\bibitem{r11}
O. Kihel and J. Larone,
\emph{On the nonnegative integer solutions to the equation $F_n\pm F_m=y^a$},
Quaest. Math. \textbf{44} (2021), no. 8, 1133--1139.
\bibitem{rnew} T. Lengyel, The order of Fibonacci and Lucas numbers, Fibonacci Quart. \textbf{33} (1995), 234–239.
\bibitem{r12}
F. Luca and V. Patel,
\emph{On perfect powers that are sums of two Fibonacci numbers},
J. Number Theory \textbf{189} (2018), 90--96.

\bibitem{r13}
E. M. Matveev,
\emph{An explicit lower bound for a homogeneous rational linear form in logarithms of algebraic numbers. II},
Izv. Ross. Akad. Nauk Ser. Mat. \textbf{64} (2000), no. 6, 125--180.
\bibitem{r14}
Z. Siar and R. Keskin,
\emph{On the Diophantine equation $F_n-F_m=2^a$},
Colloq. Math. \textbf{159} (2020), no. 1, 119--126.
\bibitem{r15}
P. Tiebekabe and I. Diouf,
\emph{On solutions of the Diophantine equation $F_{n_1}+F_{n_2}+F_{n_3}+F_{n_4}=2^a$},
J. Algebra Relat. Topics \textbf{9} (2021), no. 2, 131--148.

\bibitem{r16}
P. Tiebekabe and I. Diouf,
\emph{On the Diophantine equation $\sum_{k=1}^{5}F_{n_k}=2^a$},
Fibonacci Quart. \textbf{60} (2022), no. 5, 384--400.

\bibitem{r17}
V. Ziegler,
\emph{Sums of Fibonacci numbers that are perfect powers},
Quaest. Math. \textbf{46} (2023), no. 8, 1717--1742.

\end{thebibliography}
\end{document}